\documentclass{amsart}
\title[Towards characterising polynomiality...]{Towards characterising polynomiality of $\frac{1-q^b}{1-q^a}{n\brack m}$ and applications}
\usepackage{amssymb,amsmath,amsthm,epsfig,graphics,latexsym}
\usepackage{enumerate}
 \theoremstyle{definition}
 \newtheorem{definition}{Definition}
  \theoremstyle{plain}

  \newtheorem{theorem}    {Theorem}
  
  \newtheorem{corollary}  {Corollary}

  \newtheorem{rmk}{Remark}
\begin{document}
  \author{Mohamed El Bachraoui}
 \email{melbachraoui@uaeu.ac.ae}
  \begin{abstract}
  In this note we shall give conditions which guarantee that 
  $\frac{1-q^b}{1-q^a}{n\brack m}\in\mathbb{Z}[q]$ holds. We shall provide a full characterisation for $\frac{1-q^b}{1-q^a}{ka\brack m}\in\mathbb{Z}[q]$.
   This unifies a variety of results already known in literature.
  We shall prove new divisibility properties for the binomial coefficients and a new divisibility result 
  for a certain finite sum involving the roots of the unity.
  \end{abstract}
  \date{\textit{\today}}
  \maketitle
\section{Introduction}
\noindent
Throughout, let $\mathbb{N}$ denote the set of positive integers, let $\mathbb{N}_0= \mathbb{N}\cup\{0\}$
be the set of nonnegative integers, and let $\mathbb{Z}$ denote the set of integers. Accordingly, let
$\mathbb{Z}[q]$ denote the set of polynomials in $q$ with coefficients in $\mathbb{Z}$ and let $\mathbb{N}_0[q]$
be the set of polynomials in $q$ with coefficients in $\mathbb{N}$.
Recall that for a complex number $q$ and a complex variable $x$, the $q$-shifted factorials are given by
\[
(x;q)_0= 1,\quad (x;q)_n = \prod_{i=0}^{n-1}(1-x q^i),\quad
(x;q)_{\infty} = \lim_{n\to\infty}(x;q)_n =\prod_{i=0}^{\infty}(1-x q^i)
\]
and the $q$-binomial coefficients are given for any $m, n\in\mathbb{N}_0$ by
\[
{n\brack m} = \begin{cases}
\frac{(q;q)_n}{(q;q)_m (q;q)_{n-m}} , & \text{if\ } n\geq m\geq 0,\\
0, & \text{otherwise.}
\end{cases}
\]
Andrews~\cite{Andrews-2} introduced the function
\[
A(n,j)=\frac{1-q}{1-q^n}{n\brack j},
\]
which, for our purposes, we extend as follows.
\begin{definition}
For $a\in\mathbb{N}$ and $b, m, n \in\mathbb{N}_0$, let
\[
A(b,a;n,m) = \frac{1-q^b}{1-q^a} {n\brack m},\quad a\in\mathbb{N},\ b, m, n \in\mathbb{N}_0.
\]
We say that $A(b,a;n,m)$ is \emph{reduced} (or \emph{in reduced form}) if $a\leq n<2a$ and $0\leq m<a$.
Writing $m=u a + r$ and $n=v a + s$ with $0\leq r <a$ and $a\leq s <2a$, it is clear that the reduced form of
$A(b,a;n,m)$ is $A(b,a;s,r)$.
\end{definition}
\begin{rmk}\label{Rmk-1}
\emph{By Guo and Krattenthaler~\cite[Lemma 5.1]{Guo-Krattenthaler}, if $b\leq a$ and
$A(b,a;n,m) \in \mathbb{Z}[q]$, then $A(b,a;n,m)\in \mathbb{N}_0[q]$.}
\end{rmk}
\noindent
Slightly modifying \cite[Theorem 5]{Bachraoui-Suwaidi}, we shall show that $A(b,a;n,m)\in\mathbb{N}_0[q]$
if and only if $A(b,a;s,r)\in\mathbb{N}_0[q]$. More specifically, we have:
\begin{theorem}\label{Th-generating}
Let $a\in\mathbb{N}$ and $b, m, n \in\mathbb{N}_0$ such that $m\leq n$.
Then
\[
A(b,a;n,m) \in \mathbb{Z}[q]\
\text{if and only if\ }
A(b,a;n+la, m+ka)
 \in \mathbb{Z}[q]
\]
for all integers $k, l$ such that $0\leq m+ka\leq n+la$.
\end{theorem}
\noindent
By Theorem~\ref{Th-generating} and Remark~\ref{Rmk-1} we have:
\begin{corollary}\label{Cor-generating}
Let $a\in\mathbb{N}$ and $b, m, n \in\mathbb{N}_0$ such that $b\leq a$ and $m\leq n$.
Then
\[
A(b,a;n,m) \in \mathbb{N}_0[q]\
\text{if and only if\ }
A(b,a;n+la, m+ka)
 \in \mathbb{N}_0[q]
\]
for all integers $k, l$ such that $0\leq m+ka\leq n+la$.
\end{corollary}
\noindent
Andrews ~\cite[Theorem 2]{Andrews-2} gave the following characterisation:
\begin{equation}\label{Thm-Andrews}
A(1,n;n,m)\in \mathbb{N}_0[q]\ \text{if and only if\ } \gcd(n,m)=1.
\end{equation}
\noindent
Sun~\cite[Theorem 1.1]{Sun} proved that
\[
{an+bn \choose an} \equiv 0 \mod{\frac{bn+1}{\gcd(a,bn+1)}}.
\]
To extend this congruence, Guo and Krattenthaler~\cite[Lemma 5.2]{Guo-Krattenthaler} proved the following $q$-analogue.
\begin{equation} \label{Thm-GuoKratt-1}
A(\gcd(a,b),a+b; a+b, a)\in\mathbb{N}_0[q].
\end{equation}
\noindent
Moreover, by Guo and Krattenthaler~\cite[Theorem 3.2]{Guo-Krattenthaler} we have:
\begin{equation}\label{Thm-GuoKratt-2}
 A(\gcd(k,n),n;2n,n-k) \in\mathbb{N}_0[q]\ \text{and\ } A(k,n;2n,n-k) \in\mathbb{N}_0[q].
\end{equation}
\noindent
Notice that the functions in (\ref{Thm-Andrews}), (\ref{Thm-GuoKratt-1}), and (\ref{Thm-GuoKratt-2}) are of type
$A(b,a;n,m)$ with $a\mid n$. So, it is natural to ask for conditions guaranteeing the statement $A(b,a;na,m) \in\mathbb{N}_0[q]$ to hold.
To this end, we have the following characterisation.
\begin{theorem}\label{MainTh-1}
Let $a$, $b$, $m$, and $n$ be nonnegative integers such that $a>0$ and $na\geq m$. Then
$A(b,a;na,m) \in \mathbb{Z}[q]$ if and only if $\gcd(a,m)\mid b$.
\end{theorem}
\noindent
Combining Remark~\ref{Rmk-1} with Theorem~\ref{MainTh-1} we have the following consequence.
\begin{corollary}\label{MainCor}
Let $a$, $b$, $m$, and $n$ be nonnegative integers such that $a>0$, $b\leq a$ and $na\geq m$. Then
$A(b,a;na,m) \in \mathbb{N}_0[q]$ if and only if $\gcd(a,m)\mid b$.
\end{corollary}
\noindent
Further, Guo and Krattenthaler~\cite[Theorem 3.1]{Guo-Krattenthaler} showed that all of the functions
\begin{equation} \label{Thm-GuoKratt-3}
\begin{split}
A(1,6n-1;12n,3n),\ A(1,6n-1;12n,4n),\ A(1,30n-1;60n,6n) \\
A(1,30n-1;120n,40n),\ A(1,30n-1;120n,45n),\ A(1,66n-1;3300n,88n)
\end{split}
\end{equation}
are in $\mathbb{N}_0[q]$.
\begin{rmk}\label{Rmk-2}
\emph{To investigate the polynomiality of $A(1,a;n,m)$ we may assume by virtue of
Theorem~\ref{Th-generating} that $A(1,a;n,n-m)$ is reducible, i.e. $n=a+r$ and $n-m=a-s$ with
$0\leq r<a$ and $0\leq s<a$. In this case we have $m=r+s$ and so, we may assume that
$n=a+r$ and $n+a\geq m\geq r$.}
\end{rmk}
\noindent
Observe that the reduced forms of all of the functions listed in (\ref{Thm-GuoKratt-3}) have the form $A(1,a;a+r,m)$ with $r\leq m$. We have the following unifying argument.
\begin{theorem}\label{Unify}
Let $a\in\mathbb{N}$, let $a> r\in\mathbb{N}_0$, let $n=a+r$, and let $m\in\mathbb{N}_0$
such that $n\geq m\geq r$.
If $\gcd(a,m)=1$ and $\gcd(a,m-j)\mid n$ for all $j=1,\ldots,r$, then
$A(1,a;n,m)\in \mathbb{N}_0[q]$.
\end{theorem}
\noindent
For instance, applying Theorem~\ref{Unify} to $a=6n-1$, $r=2$, and $m=3n$ gives that $A(1,6n-1;12n,3n) \in \mathbb{N}_0[q]$
and applying Theorem~\ref{Unify} to $a=30n-1$, $r=4$, and $m=45n$ gives that $A(1,30n-1;120n,45n) \in \mathbb{N}_0[q]$.
One can check the polynomiality of the other functions listed in (\ref{Thm-GuoKratt-3}) in a similar way.

\noindent
An important application of the function $A(b,a;n,m)$ is the fact that whenever it is a polynomial in
$\mathbb{Z}[q]$ and $\gcd(a,b)=1$, then $a\mid{n\choose m}$.
Our next result deals with divisibility properties for the binomial coefficients.
\begin{theorem}\label{binom-div}
If $a$ and $n$ are nonnegative integers such that $a\geq 3$, then
\[ \emph{(a)}\quad \bigl( (a-1)n+1\bigr) \Big\vert
\gcd\left( {(a-1)^2 n-1\choose (a-1)n},  {a(a-1) n\choose 2(a-1)n+1} \right),\]
\[ \emph{(b)}\quad \bigl( (a-1)n-1\bigr) \Big\vert
\gcd \left({(a-1)^2 n-1\choose (a-1)n-2}, {a(a-1) n-2\choose 2(a-1)n-3} \right) .  \]
\end{theorem}
\noindent
Finally, by a result of
Gould~\cite{Gould} we have for any nonnegative integers $N$ and $M<n$
\begin{equation}\label{Gould-sum}
\sum_{j\geq0}{N+mn\choose M+jn} = \frac{1}{n}\sum_{j=1}^{n} w^{-jM}(1+ w^j)^{N+mn},
\end{equation}
where $w=e^{2\pi i/n}$ is a primitive $n$th root of unity.
In particular, this implies  that
\[
n\Big\vert \sum_{j=1}^{n} w^{-jM}(1+ w^j)^{N+mn}.
\]
We have the following generalisation.
\begin{theorem}\label{application-root-unity}
If $A(1,n;N,M)\in\mathbb{Z}[q]$ , then for any nonnegative integer $m$ we have
\[
  n^2 \Big\vert \sum_{j=1}^{n} w^{-jM}(1+ w^j)^{N+mn},
\]
where $w=e^{2\pi i/n}$ is a primitive $n$th root of unity.
\end{theorem}
\section{Proof of Theorem \ref{Th-generating}}
The implication from the right to the left is clear. Assume now that
\[
A(b,a;n,m) \in \mathbb{Z}[q].
\]
By the well-known identity
\[
q^M -1 = \prod_{d\mid M}\Phi_d(q),
\]
where $\Phi_d(q)$ is the $d$-th cyclotomic polynomial in $q$, we obtain
\[
A(b,a;n,m) = \prod_{d=2}^{n} \Phi_d(q)^{e_d},
\]
where
\[
e_d = \chi(d\mid b)-\chi(d\mid a) + \left\lfloor\frac{n}{d}\right\rfloor - \left\lfloor\frac{m}{d}\right\rfloor - \left\lfloor\frac{n-m}{d}\right\rfloor,
 \]
with $\chi(S)=1$ if $S$ is true and $\chi(S)=0$ if $S$ is false. As $A(b,a;n,m)\in\mathbb{Z}[q]$ and $\Phi_d(q)$ is irreducible for any $d$ we must have
$e_d \geq 0$ for all $d=2,\ldots,n$. As to
$A(b,a;n+la, m+ka)$,
we have
\[
A(b,a;n+la, m+ka) = \prod_{d=2}^{n+la} \Phi_d(q)^{e_d},
\]
where
\[
e_d = \chi(d\mid b)-\chi(d\mid a) + \left\lfloor\frac{n+la}{d}\right\rfloor - \left\lfloor\frac{m+ka}{d}\right\rfloor
 - \left\lfloor\frac{n-m+(l-k)a}{d}\right\rfloor,
 \]
Then clearly $e_d \geq 0$ unless $d\mid a$. But if $d\mid a$, then
\[
\left\lfloor\frac{n+la}{d}\right\rfloor - \left\lfloor\frac{m+ka}{d}\right\rfloor
 - \left\lfloor\frac{n-m+(l-k)a}{d}\right\rfloor
=
\left\lfloor\frac{n}{d}\right\rfloor - \left\lfloor\frac{m}{d}\right\rfloor
 - \left\lfloor\frac{n-m}{d}\right\rfloor
\]
and therefore $e_d \geq 0$ by assumption, implying that $A(b,a;n+la, m+ka)$ is a polynomial in $q$.
\section{Proof of Theorem~\ref{MainTh-1}}
Suppose that $\gcd(a,m)=g\nmid b$ and that $A(b,a;na,m)\in\mathbb{Z}[q]$. Then clearly
$A(b,g;na,m)\in\mathbb{Z}[q]$ and so, by Theorem~\ref{Th-generating} we have
\[
A(b,g;na,0)= \frac{1-q^b}{1-q^g}\in\mathbb{Z}[q],
\]
which is impossible as $g\nmid b$.
Assume now that $\gcd(a,m)\mid b$. Then just as before, we have
\[
A(b,a;na,m) = \prod_{d=2}^{na} \Phi_d(q)^{e_d},
\]
where
\[
e_d = \chi(d\mid b)-\chi(d\mid a) + \left\lfloor\frac{na}{d}\right\rfloor - \left\lfloor\frac{m}{d}\right\rfloor - \left\lfloor\frac{na-m}{d}\right\rfloor.
 \]
 Then $e_d \geq 0$ unless $d\mid a$. But if $d\mid a$, then
 \begin{equation}\label{floor}
 e_d = \chi(d\mid b)-1-\left(\left\lfloor\frac{m}{d}\right\rfloor + \left\lfloor\frac{-m}{d}\right\rfloor\right).
 \end{equation}
\noindent
 {\bf Case 1:} $d\mid m$. Then $d\mid\gcd(a,m)$ and so also $d\mid b$. From these facts and
 the identity (\ref{floor}) we conclude that $e_d=0$.

\noindent
 {\bf Case 2:} $d\nmid m$. Then $\left\lfloor m/d \right\rfloor + \left\lfloor -m/d \right\rfloor = -1$ and so,
 $e_d= \chi(d\mid b) -1+1 \geq 0$. This completes the proof.
 \section{Proof of Theorem~\ref{Unify}}
Proceeding as before, we have
\[
A(1,a;n,m) = \prod_{d=2}^{n} \Phi_d(q)^{e_d},
\]
with
\[
e_d = -\chi(d\mid a) + \left\lfloor\frac{n}{d}\right\rfloor - \left\lfloor\frac{m}{d}\right\rfloor -
\left\lfloor\frac{n-m}{d}\right\rfloor.
\]
 Then $e_d \geq 0$ unless $d\mid a$. Let  $2\leq d\mid a$. Suppose that there is some
 $j=1,\ldots,r$ such that $d\mid \gcd(a,m-j)$. Then $d\mid n$ but $d\nmid m$ and we get
\[
e_d = -1 - \left(\left\lfloor\frac{m}{d}\right\rfloor+\left\lfloor\frac{-m}{d}\right\rfloor\right) = -1- (-1) \geq 0.
\]
Suppose now that
$d\nmid m-1, \ldots, d\nmid m-r$. Then
\[
\left\lfloor\frac{m}{d}\right\rfloor + \left\lfloor\frac{r-m}{d}\right\rfloor = \left\lfloor\frac{m}{d}\right\rfloor + \left\lfloor\frac{-m}{d}\right\rfloor = -1
\]
and so,
\[
\begin{split}
e_d &= -1 + \left\lfloor\frac{a+r}{d}\right\rfloor -
\left( \left\lfloor\frac{m}{d}\right\rfloor + \left\lfloor\frac{a+r-m}{d}\right\rfloor \right) \\
&= -1 + \left\lfloor\frac{r}{d}\right\rfloor -
\left( \left\lfloor\frac{m}{d}\right\rfloor + \left\lfloor\frac{r-m}{d}\right\rfloor \right)  \\
&= -1 + \left\lfloor\frac{r}{d}\right\rfloor - (-1) \\
& \geq 0,
\end{split}
\]
implying that $A(1,a;n,m) \in\mathbb{Z}[q]$. The fact that $A(1,a;a+r,m) \in\mathbb{N}_0[q]$ is a consequence of Remark~\ref{Rmk-1}.
\section{Proof of Theorem~\ref{binom-div}}
(a)\ Let $a\geq 3$ and $n$ be nonnegative integers.
From the evident fact
\[
A(an+1,n+1;an,n) = {an+1\brack n+1} \in\mathbb{Z}[q]
\]
and Theorem~\ref{Th-generating} we get
\[ A(an+1,n+1;an-n-1,n)\in\mathbb{Z}[q] \ \text{and\ } A(an+1,n+1;an,n+n+1) \in\mathbb{Z}[q],
\]
from which we find
\[
(n+1) \Big\vert (an+1)\gcd\left( {(a-1)n-1 \choose n}, {an \choose 2n+1} \right).
\]
Letting $n:= (a-1)m$ we have that $\gcd(n+1, an+1)=1$ and so the previous divisibility implies
\[
\bigl( (a-1)m+1 \bigr) \Big\vert \gcd\left( {(a-1)^2 m-1 \choose (a-1) m-1}, {a(a-1)m \choose 2(a-1)m+1} \right),
\]
as desired.

(b)\ Follows similarly by applying Theorem~\ref{Th-generating} to the fact
\[
A(an-1,n-1;an-2,n-2) = {an-1\brack n-1} \in\mathbb{Z}[q].
\]
\section{Proof of Theorem~\ref{application-root-unity}}
Suppose that $A(1,n;N,M)\in\mathbb{Z}[q]$. Then
by virtue of Theorem~\ref{Th-generating} we have
\[
A(1,n;N+mn,M+jn) \in \mathbb{Z}[q]
\]
for all nonnegative integers $j$ such that $M+jn\leq N+mn$.
It follows with the help of Gould's identity~(\ref{Gould-sum})
\[
n\Big\vert \sum_{j\geq0} {N+mn\choose N+jn} =
\frac{1}{n} \sum_{j=1}^{n} w^{-jM}(1+w^{j})^{N+mn},
\]
from which the desired divisibility follows.

%
\end{document}